\documentclass[pre,nofootinbib,twocolumn]{revtex4}
\usepackage[dvips]{hyperref}
\usepackage{amsmath,amssymb}
\usepackage[dvips]{graphics}
\DeclareMathOperator*{\argmin}{argmin} % Thanks, Wikipedia!
\newcommand{\censoredlike}{\ensuremath{\ell_{\mathrm{C}}}}
\newcommand{\censoredtheta}{\ensuremath{\hat{\theta}_{\mathrm{C}}}}
\newcommand{\censoredsigma}{\ensuremath{\hat{\sigma}_{\mathrm{C}}}}

\begin{document}
\title{Maximum Likelihood Estimation for $q$-Exponential (Tsallis) Distributions}
\author{Cosma Rohilla Shalizi}
\affiliation{Statistics Department, Carnegie Mellon University}
\email{cshalizi@cmu.edu}
\date{Begun 28 December 2006, last updated 31 January 2007}
\begin{abstract}
  This expository note describes how to apply the method of maximum likelihood
  to estimate the parameters of the ``$q$-exponential'' distributions
  introduced by Tsallis and collaborators.  It also describes the relationship
  of these distributions to the classical Pareto distributions.
\end{abstract}
\maketitle

In a series of papers beginning with \cite{Tsallis-first-paper}, Constantino
Tsallis and collaborators introduced what have come to be called
$q$-exponential probability distributions.  These can be defined through their
``complementary'' (``upper'', ``upper cumulative'') distribution functions,
also called ``survival'' functions:
\begin{equation}
\label{eqn:q-kappa-survival-function}
P_{q,\kappa}(X \geq x) = {\left(1-\frac{(1-q)x}{\kappa}\right)}^{1/(1-q)}
\end{equation}
Tsallis {\em et al.}\ proposed these distributions to handle
statistical-mechanical systems with long-range interactions, necessitating (it
is claimed) a non-extensive generalization of the ordinary Gibbs-Shannon
entropy.  Following Jaynes's procedure of maximizing an entropy subject to
constraints on expectation values \cite{Jaynes-essays}, they got the
$q$-exponential distributions, in which $\kappa$ enforces the constraints, and
$q$ measures the departure from extensivity, Boltzmann-Gibbs statistics being
recovered as $q \rightarrow 1$.

Tsallis's ideas about non-extensive entropy and its possible applications, in
and out of statistical mechanics, have attracted intense (not to say
``extensive'') interest in physics; the bibliography at
\url{http://tsallis.cat.cbpf.br/biblio.htm} has over 2000 entries.  They are
also quite controversial (see, e.g., Refs.\
\cite{LaCourt-Schieve-contra-Tsallis,Zanette-Montemurro-temperature,%
  Zanette-Montemurro-nonthermodynamic,Bouchet-Dauxois-Ruffo-controversy,%
  Lavenda-Dunning-Davies-additive-entropies-of-degree-q,%
  Nauenberg-critique-of-q-entropy}, the replies by Tsallis and others, and in
some cases the replies to the replies).  Whether or not the critics are
correct, however, $q$-exponentials are still valid probability distributions,
and can usefully describe some empirical phenomena.  To this end, in a recent
paper Douglas R. White {\em et al.} pose the problem of estimating the
parameters $q$ and $\kappa$ from data by the method of maximum likelihood
\cite{White-Kejzar-Tambayong-oscillatory-city-sizes}.  This note solves that
problem.

I first reparameterize Eq.\ \ref{eqn:q-kappa-survival-function} to simplify
estimation and emphasize links to Pareto distributions.  I then rehearse the
math of finding the maximum likelihood estimator (MLE) for the $q$-exponential
distribution, discussing its accuracy and precision, and adjustments for data
in which samples below a fixed threshold are all dropped (``censoring'').  I
compare maximum-likelihood estimates to those found by the current practice of
curve-fitting; the latter are inferior.  Finally, I discuss testing the
assumption that the data are $q$-exponentially distributed.  Code implementing
the MLE for $q$-exponentials is available at
\url{http://bactra.org/research/tsallis-MLE/}, written in \texttt{R}, a free,
open-source programming language for statistical computing
(\url{http://www.r-project.org/}).  This code also calculates probabilities and
quantiles, generates random numbers, etc.

\paragraph{Reparameterization as Generalized Pareto Distributions}
\label{sec:q-exp-are-GPDs}

While it is possible to find the MLE for $q$-exponentials in the form given in
Eq.\ \ref{eqn:q-kappa-survival-function}, %\footnote{And I have.},
the algebra is needlessly messy.  It is simpler to reparameterize, and change
back to the original parameter system at the end, if desired.  (Under a 1-1
change of parameters, an MLE for the old parameters must, under the
transformation, be an MLE for the new parameters, and vice versa.)  Thus,
define the new parameters $\theta \equiv -\frac{1}{1-q}$ and $\sigma \equiv
\theta \kappa$, from which the original parameters can be recovered:
\begin{equation}
\label{eqn:back-to-q-kappa} q  = 1 + \frac{1}{\theta}, ~ \kappa  = \frac{\sigma}{\theta}
\end{equation}
In the new parameter system, the survival function becomes
\begin{equation}
\label{eqn:survival-function-in-sane-parameters}
P_{\theta,\sigma}(X \geq x) = {\left(1+x/\sigma\right)}^{-\theta}
\end{equation}
Hence the probability density is
\begin{equation}
p_{\theta,\sigma}(x) % & \equiv & -\frac{\partial P_{\theta,\sigma}(X \geq x)}{\partial x}\\
 = \frac{\theta}{\sigma}{\left(1+x/\sigma\right)}^{-\theta -1}
\label{eqn:q-exp-density}
\end{equation}
The code mentioned above uses both parameterizations.

$Y$ has a Pareto distribution with scaling exponent $\alpha$ and cut-off $y_0$
if $p(y) = 0$ when $y < y_0$, and otherwise $p(y) \propto
{(y/y_0)}^{-\alpha-1}$.  Hence when $X$ has a $q$-exponential distribution,
$1+x/\sigma$ has a Pareto distribution with cut-off 1 and scaling exponent
$\theta$.  Following the classification given in Arnold's monograph on Pareto
distributions \cite{Arnold-on-Pareto-distributions}, this is an instance of a
``type II generalized Pareto'', often used in operations research on failure
times and other reliability problems, the standard form of which is $P(X \geq
x) = {[1 + (x-\mu)/\sigma]}^{-\alpha}$.  The $q$-exponentials come from taking
$\mu = 0$ and $\alpha = \theta$; the ordinary Pareto distribution is recovered
by taking $\sigma = x_0$ and $\mu=\sigma$.  According to Ref.\ \cite[pp.\
13--14, 208--210]{Arnold-on-Pareto-distributions}, the type II generalized
Pareto was introduced in Refs.\
\cite{Maguire-Pearson-Wynn,Silcock-labor-turnover,Harris-Pareto-queues}, and
the latter two also derived the MLE.\footnote{Arnold \cite[p.\
  48]{Arnold-on-Pareto-distributions} observes that a mixture of exponentials
  can produce a type II generalized Pareto.  If the distribution of $X-\mu$,
  given $Z$, is an exponential with mean $\sigma/Z$, and $Z$ has a
  $\Gamma(\alpha,1)$ distribution, then $X$ has a type II generalized Pareto
  distribution with parameters $\mu$, $\sigma$ and $\alpha$.  He assigns
  priority for this result to \cite{Maguire-Pearson-Wynn}.  It would appear to
  be equivalent to C. Beck's ``superstatistics'' approach to Tsallis statistics
  (reviewed in \cite{Beck-superstatistics-recent-developments}).}  The
calculations below are a special case of their results, except for the
treatment of censoring, which may be new.

\paragraph{Derivation of the MLE for $q$-Exponentials}
\label{sec:MLE}

Under the $q$-exponential model with parameters $\theta,\sigma$, the
log-probability density of a sequence of independent, identically-distributed
samples $X_1 = x_1$, $X_2 = x_2$, \ldots $X_n = x_n$, for short $X_1^n =
x_1^n$, is
\begin{eqnarray}
\log{p_{\theta,\sigma}(x_1^n)} & = &  -n\log{\sigma} + n\log{\theta} \\
\nonumber & & - (\theta+1)\sum_{i=1}^{n}{\log{1+x_i/\sigma}}\\
& \equiv & \ell(\theta,\sigma) ~,
\end{eqnarray}
the log-likelihood of the parameter combination $\theta,\sigma$.

To find the MLEs, take the first derivatives of the log-likelihood with respect
to the parameters and set them equal to zero.  First, the shape parameter
$\theta$:
\begin{eqnarray}
\frac{\partial \ell}{\partial \theta} % & = & 
% \frac{n}{\theta} - \left(\frac{\partial (\theta +1)}{\partial \theta}\right)\left(\sum_{i=1}^{n}{\log{1+x_i/\sigma}}\right)\\
%\nonumber & &  - (\theta+1)\left(\sum_{i=1}^{n}{\frac{1}{1+x_i/\sigma}\frac{\partial (1+x_i/\sigma)}{\partial \theta}}\right)\\
& = & \frac{n}{\theta} - \sum_{i=1}^{n}{\log{1+x_i/\sigma}}\\
% 0 & = & \frac{n}{\hat{\theta}} - \sum_{i=1}^{n}{\log{1+x_i/\sigma}}\\
\hat{\theta} & = & n{\left[\sum_{i=1}^{n}{\log{1+x_i/\sigma}}\right]}^{-1}
\label{eqn:theta-hat-in-terms-of-sigma}
\end{eqnarray}
Similarly for the scale parameter $\sigma$:
\begin{eqnarray}
\frac{\partial \ell}{\partial \sigma} % & = & -\frac{n}{\sigma} - (\theta +1)\sum_{i=1}^{n}{\frac{1}{1 + x_i/\sigma} \frac{\partial (1+x_i/\sigma)}{\partial \sigma}}\\
& = & -\frac{n}{\sigma} + \frac{\theta+1}{\sigma^2}\sum_{i=1}^{n}{\frac{x_i}{1 + x_i/\sigma}}\\
% 0 & = & n\hat{\sigma} - (\theta+1)\sum_{i=1}^{n}{\frac{x_i}{1 + x_i/\hat{\sigma}}}\\
\hat{\sigma} & = & \frac{\theta+1}{n}\sum_{i=1}^{n}{\frac{x_i}{1 + x_i/\hat{\sigma}}}
\label{eqn:sigma-hat-in-terms-of-theta}
\end{eqnarray}

Eqs.\ \ref{eqn:theta-hat-in-terms-of-sigma} and
\ref{eqn:sigma-hat-in-terms-of-theta} give the MLEs for $\theta$ and $\sigma$,
respectively, if the other parameter is known.  The former gives the value of
$\hat{\theta}$ explicitly\footnote{Cf.\ the well-known MLE for the scaling
  exponent in a Pareto distribution
  \cite{Muniruzzaman-on-Pareto,Arnold-on-Pareto-distributions,%
    MEJN-on-power-laws}, $\hat{\alpha} =
  n/\left[\sum_{i=1}^{n}{\log{x/x_0}}\right]$.}, while the latter does so
implicitly, through the solution of an equation.  Implicitly-defined MLEs like
this occur in several generalizations of the exponential distribution, such as
the ones known to physicists as ``stretched exponentials'' and to statisticians
as ``Weibull distributions'' (after the physicist who introduced them)
\cite[ch.\ 20]{Johnson-Kotz-continuous-univariate-1}.  The lack of a closed
form is only a small annoyance, since such equations can generally be rapidly
solved numerically, to a precision much smaller than the uncertainty inherent
in the data.

If neither $\theta$ nor $\sigma$ is known (i.e., neither $q$ nor $\kappa$),
then the simultaneous solution of Eqs.\ \ref{eqn:theta-hat-in-terms-of-sigma}
and \ref{eqn:sigma-hat-in-terms-of-theta} gives the joint maximum likelihood
estimator.  Substituting the former equation into the latter gives a single
equation in $\hat{\sigma}$ and the data:
\begin{equation}
\hat{\sigma}  =  \frac{1}{n}\left(1 + n{\left[\sum_{i=1}^{n}{\log{1+x_i/\hat{\sigma}}}\right]}^{-1}\right)\sum_{i=1}^{n}{\frac{x_i}{1+x_i/\hat{\sigma}}}
\label{eqn:sigma-hat-on-its-own}
\end{equation}
This does not seem to simplify, but, again, can be solved numerically.  (Eq.\
\ref{eqn:sigma-hat-on-its-own} is transcendental, whereas Eq.\
\ref{eqn:sigma-hat-in-terms-of-theta} is rational, but no worse than the
equation for the MLE of the Weibull distribution, which also contains a sum of
logarithms, etc.)  Substituting the solution into Eq.\
\ref{eqn:theta-hat-in-terms-of-sigma} gives $\hat{\theta}$, and then Eq.\
\ref{eqn:back-to-q-kappa} give $\hat{q}$, $\hat{\kappa}$.

\paragraph{Accuracy and Precision of the MLE}
\label{sec:accuracy}

An estimator $\hat{\psi}(X_1^n)$ of a parameter $\psi$ of a statistical
distribution is {\em consistent} when $\hat{\psi}$ converges in probability to
$\psi$, i.e., for any $\epsilon > 0$ and any $\delta > 0$, for sufficiently
large $n$, $P\left(\left\|\hat{\psi}(X_1^n) - \psi\right\| \geq \epsilon\right)
\leq \delta$.  In other words, a consistent estimator is ``probably
($1-\delta$) approximately ($\epsilon$) correct'', for arbitrarily small
$\delta$ and $\epsilon$.  Under quite general conditions, met here,
maximum likelihood estimators are consistent \cite{Pitman-basic-theory}.

Consistency alone is not enough to calculate standard errors or confidence
regions.  However, under conditions only mildly more restrictive than those
needed for consistency, MLEs are {\em asymptotically normal} and {\em
  unbiased}.  That is, $\hat{\psi}(X_1^n) - \psi$ has, for large $n$, a
multidimensional Gaussian distribution with mean zero and covariance matrix
$(1/n) I^{-1}(\psi)$, where $I(\psi)$ is the {\em Fisher information matrix},
\begin{equation}
\label{eqn:fisher-info-defined} I_{ij}(\psi) \equiv -\int{\frac{\partial^2 \log{p_{\psi}(x)}}{\partial \psi_{i} \partial \psi_{j}} p_{\psi}(x) dx}
\end{equation}
By the famous Cram\'er-Rao inequality \cite{Cramer}, any consistent unbiased
estimator has a covariance at least equal to $I^{-1}(\psi)$; the MLE is {\em
  asymptotically efficient} because it attains this bound.  Since the true
value of $\psi$ is unknown, $I(\psi)$ cannot give us standard errors or
confidence regions, but $I(\hat{\psi})$ is a consistent estimator of $I(\psi)$,
and can be used for those purposes.  Another consistent estimator of the Fisher
information is the {\em observed information matrix}, $J_{ij}(\psi) \equiv -
n^{-1}\partial^2 \ell(\psi)/\partial \psi_{i} \partial \psi_{j}$,
%\begin{equation}
%J_{ij}(\psi) = -\frac{1}{n}\frac{\partial \ell(\psi)}{\partial \psi_{i} \partial \psi_{j}} 
%\end{equation}
and $J(\hat{\psi})$ also gives asymptotically-correct error estimates.  Ref.\
\cite{Barndorff-Nielsen-and-Cox-inference-and-asymptotics} treats these
standard results in detail.

For $q$-exponential distributions, it is easy to verify that the standard
conditions for the asymptotic normality of the MLE hold.  In the
$\theta,\sigma$ parameterization, simple but lengthy calculus yields
\begin{equation}
\label{eqn:fisher-information} I(\theta,\sigma) = \left[ \begin{array}{cc} \frac{1}{\theta^2} & -\frac{1}{(\theta+1)\sigma} \\ -\frac{1}{(\theta+1)\sigma} & \frac{\theta}{\sigma^2(\theta+2)} \end{array} \right]
\end{equation}
Either $I(\hat{\theta},\hat{\sigma})$ or the observed information matrix could
be used to find standard errors and Gaussian confidence regions.  Propagation
of errors can then carry these to estimates on $q$ and $\kappa$.

For small samples, asymptotic approximations should be avoided in favor of {\em
  parametric bootstrapping} \cite[sec.\ 9.11]{Wasserman-all-of-stats}.  Having
obtained an estimate $\hat{\psi} = \hat{\psi}(x_1^n)$, make up a ``bootstrap''
sample of random numbers $Y_1, Y_2, ... Y_n$ with the density $p_{\hat{\psi}}$,
and calculate $\hat{\psi}(Y_1^n)$. The distribution of $\hat{\psi}(Y_1^n) -
\hat{\psi}$ is approximately the same as that of $\hat{\psi}(X_1^n) - \psi$, so
by taking many bootstrap samples one can estimate standard errors and
confidence regions, without making Gaussian approximations.  (For more on
bootstrapping, see, e.g., \cite[ch.\ 8]{Wasserman-all-of-stats}.)  The code
mentioned above finds bootstrapped biases, standard errors and confidence
intervals.

\paragraph{Censored Data}
\label{sec:censored}

In many applications, only measurements exceeding some known lower threshold
$x_0$ are available, i.e., only values of $X \geq x_0$ become data.  Parameters
estimation from such {\em left-censored} data must take account of the
threshold.  Specifically, rather than maximizing the unconditional likelihood,
$\ell(\theta,\sigma)$, one should maximize the likelihood conditional on being
in the right tail, $\censoredlike(\theta,\sigma,x_0)$.  It is easily shown that
the censored density is 0 when $x < x_0$, and otherwise
\begin{equation}
p_{\theta,\sigma,x_0}(x) = {\left(1+x_0/\sigma\right)}^{\theta} p_{\theta,\sigma}(x)
\end{equation}
$p_{\theta,\sigma}(x)$ being given by Eq.\ \ref{eqn:q-exp-density}.  The
censored likelihood thus equals $\ell(\theta,\sigma)$ plus a term involving
only $\theta$, $\sigma$ and $x_0$:
\begin{equation}
\label{eqn:censored-likelihood} \censoredlike(\theta,\sigma,x_0) = \ell(\theta,\sigma) + n\theta\log{1+x_0/\sigma}
\end{equation}
The likelihood estimating equations become
\begin{eqnarray}
\label{eqn:censoredtheta} \censoredtheta & = & n{\left[\sum_{i=1}^{n}{\log{\frac{1+x_i/\sigma}{1+x_0/\sigma}}} \right]}^{-1}\\
\label{eqn:censoredsigma} \censoredsigma & = &  - \theta\frac{x_0}{1+x_0/\censoredsigma} + \frac{\theta+1}{n}\sum_{i=1}^{n}{\frac{x_i}{1+x_i/\censoredsigma}}
\end{eqnarray}
Eqs.\ \ref{eqn:censoredtheta} and \ref{eqn:censoredsigma} reduce to Eqs.\
\ref{eqn:theta-hat-in-terms-of-sigma} and \ref{eqn:sigma-hat-in-terms-of-theta}
when $x_0 = 0$ (no censoring), and can be solved in the same way.  The MLE
remains consistent, and asymptotically normal and efficient.  The Fisher
information matrix, after an even longer calculation, ends up being $I(\theta,\sigma+x_0)$; explicitly,
\begin{equation}
I_{\mathrm{C}}(\theta,\sigma,x_0) = \left[ \begin{array}{cc} \frac{1}{\theta^2} & -\frac{1}{(\theta+1)(\sigma+x_0)} \\ -\frac{1}{(\theta+1)(\sigma+x_0)} & \frac{\theta}{{(\sigma+x_0)}^2(\theta+2)} \end{array} \right]
\end{equation}
Bootstrapping, however, is even more strongly recommended than with uncensored
data.  Simulated values should be drawn from the tail only.

\paragraph{Comparison to Curve-Fitting}
\label{sec:least-squares}

Hitherto, attempts to estimate the parameters of $q$-exponential distributions
have been based on curve-fitting.  (Ref.\
\cite{White-Kejzar-Tsallis-Farmer-White} is an unusually careful example.)
Taking the log of both sides of Eq.\
\ref{eqn:survival-function-in-sane-parameters},
\begin{equation}
\log{P_{\theta,\sigma}(X\geq x)} = -\theta \log{\left(1+x/\sigma\right)}
\end{equation}
Write $S_n(x)$ for the empirical distribution function, i.e., the fraction of
points in $x_1, x_2, \ldots x_n$ which are $\geq x$.  Then we expect that, at
least for large sample sizes $n$,
\begin{eqnarray}
\label{eqn:empirical-survival-function}
\log{S_n(x_i)} & \approx & -\theta \log{\left(1+x_i/\sigma\right)}
\end{eqnarray}
for all $x_i$.  A least-squares approach to estimating the parameters minimizes
the squared difference between the two sides of Eq.\
\ref{eqn:empirical-survival-function}, summed over all $x_i$, i.e.,
\begin{equation}
% \lefteqn{(\tilde{\theta},\tilde{\sigma}) = } & &\\
(\tilde{\theta},\tilde{\sigma})  \equiv  \argmin_{\theta,\sigma}{\sum_{i=1}^{n}{{\left(\log{S_n(x_i)} + \theta \log{\left(1+\frac{x_i}{\sigma}\right)}\right)}^2}}
\end{equation}
It is possible to show that this estimator is consistent.

\begin{figure}[t]
\resizebox{\columnwidth}{!}{\includegraphics{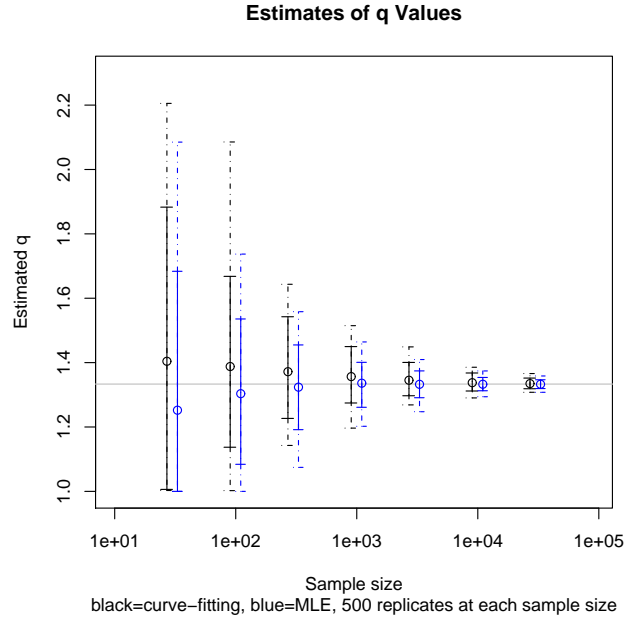}}
\caption{(Color online) Comparison of estimates of $q$ using the MLE and
  curve-fitting.  All data generated using $q = 4/3, \kappa=200/3$
  ($\theta=3,\sigma=200$), with varying sample sizes, and 500 independent
  replications are each sample size. Curve-fitting estimates are plotted in
  black, displaced slightly to the left, and MLEs in blue, displaced to the
  right.  Solid bars show the 5th and 95th percentiles of sample estimates,
  circles the median estimate, and dashed lines the sample extrema.  Note that
  the MLE is always less biased and more precise.}
\label{fig:curve-fitting-is-bad}
\end{figure}

The analogous procedure for Pareto distributions was the one originally used by
Pareto in the 1890s \cite{Arnold-on-Pareto-distributions}, and still widespread
in physics.  For Pareto distributions, however, statisticians have known since
the 1950s that such estimation-by-regression is much more biased, and much less
precise, than the maximum likelihood estimator
\cite{Muniruzzaman-on-Pareto,Arnold-on-Pareto-distributions}.  (In particular,
the standard errors are much larger than blind use of the ordinary regression
formulas suggest.)  The same is true of the least-squares estimate of
$q$-exponentials (Fig.\ \ref{fig:curve-fitting-is-bad}).  Fitting curves to
binned estimates of the probability density, rather than to the cumulative
distribution, is even less accurate.  Neither approach should be used.

\paragraph{Validation}
\label{sec:validation}

All the claims of consistency, efficiency, etc., made above assume that the
data really do come from a $q$-exponential distribution.  In statistical
terminology, the assumption is that the $q$-exponential model is correctly {\em
  specified}, as opposed to being {\em mis-specified}.  In applications to
empirical data, it is crucial to check this assumption.  Rigorous
mis-specification tests are too complicated to go into here
\cite{White-specification-analysis,Spanos-textbook}, but some remarks are in
order.

The most common test of specification in the literature on Tsallis statistics
is to look at the fraction, $R^2$, of the variance in $\log{S_n}$ accounted for
by the fitted distributional curve.  Unfortunately, this popularity is not
based on any reliability; it is easy to construct examples where $1+x/\sigma$
has, say, a log-normal distribution, but $R^2$ is always close to 1.  Rather
than looking at $R^2$, one should either test $q$-exponentials against
alternative distributions such as the Pareto, the log-normal, etc., or do
general goodness-of-fit tests, adjusting for the way parameters are estimated
from the data \cite[ch.\ 10]{Wasserman-all-of-stats}.  The latter must be
interpreted with caution: failing a goodness-of-fit test provides strong
evidence against a model, but passing one may give only very weak evidence in
its favor, depending on the severity of the test \cite{Mayo-Cox-frequentist}.

Two heuristic checks for mis-specification deserve mention.  One compares the
parametric bootstrap, described above, with a {\em non-parametric bootstrap},
in which the values $Y_1, \ldots Y_n$ come from resampling the data $x_1,
\ldots x_n$ with replacement, not from the fitted distribution.  If parametric
and non-parametric bootstrap estimates of bias, standard error, etc., differ
substantially, this is a sign that the model is mis-specified.  Similarly, if
the expected Fisher information at the MLE, $I(\hat{\theta},\hat{\sigma})$ is
very different from the observed information, $J(\hat{\theta},\hat{\sigma})$,
this again suggests the statistical model poorly describes the data-generating
process.  The comparison of information matrices can be turned into a formal
test for mis-specification \cite{White-specification-analysis}.

\paragraph*{Conclusion}

Tsallis $q$-exponentials are legitimate possible models of heavy-tailed data.
Under other names, they have been so used in operations research and statistics
for half a century, without any entropic origin story.  To model data with
$q$-exponentials, their parameters must be estimated accurately.  The
estimators currently used by physicists are inferior to the MLE, which is
asymptotically efficient.  If physicists want to describe data with
$q$-exponentials, they should stop fitting curves and start maximizing
likelihoods.  Whether using Tsallis statistics is a good idea in the first
place is another matter, beyond the scope of this note.

{\em Acknowledgments.}  Thanks to Doug White for valuable discussions, sharing
manuscripts, and pressing me to work through this problem and write it up.

\bibliography{locusts}

\end{document}